\documentclass[12pt]{article}
\usepackage{amssymb}
\usepackage{amsthm}\usepackage{epsf,epsfig,amsfonts}
\usepackage{amsmath}
\usepackage{graphicx}

\usepackage{a4wide}

\newtheorem{corollary*}{Corollary}

\newcommand{\be}{\begin{equation}}
\newcommand{\ee}{\end{equation}}

\newcommand{\weg}[1]{}

\newcommand{\C}[1]{\stackrel{#1}{c}}
\usepackage{epsf,epsfig,amsfonts,a4wide}

\usepackage{amsfonts}
\usepackage{amsmath,amssymb,amsthm}
\usepackage{url}
\usepackage{amsmath,amssymb,amsthm}

\usepackage{amsfonts}
\usepackage{amsmath,amssymb,amsthm}
\usepackage{url}

\newtheorem{Th}{Theorem}

\newtheorem{Lemma}{Lemma}
\newtheorem{Cor}{Corollary}

\theoremstyle{remark}

\newtheorem{Rem}{Remark}

\newcommand{\eqdef}{\stackrel{{\rm def}}{=}}

\newcommand{\const}{\mbox{\rm const}}
\title{Complete Einstein metrics are geodesically rigid}
\date{} \author{ Volodymyr  Kiosak, Vladimir  S. Matveev\thanks{ Institute of Mathematics, FSU Jena, 07737 Jena Germany,  matveev@minet.uni-jena.de}}
\begin{document}
\maketitle
\begin{abstract}
We  prove that every  complete Einstein (Riemannian or  pseudo-Riemannian)  metric $g$ is
 geodesically rigid: if any other complete metric $\bar g$  has the same 
 (unparametrized) geodesics with $g$, then the Levi-Civita connections of $g$ and $\bar g$ coincide.
\end{abstract}
\section{Introduction}
\subsection{Definitions and results}
\label{results} 

Let $(M^n,g)$ be a connected Riemannian or pseudo-Riemannian manifold of dimension $n\ge 3$.   We say that a metric $\bar g$ on $M^n$ is \emph{geodesically equivalent} to $g$, if every geodesic of $g$ is a (reparametrized) geodesic of $\bar g$.  We say that they are \emph{affine equivalent},  if their Levi-Civita connections coincide. We say that $g$ is \emph{Einstein}, 
if $R_{ij}= \tfrac{R}{n} \cdot g_{ij},$ where $R_{ij}$ is the Ricci tensor of the metric $g$, and $R:= R_{ij}g^{ij}$ is the scalar curvature.   Our main result is 

\begin{Th}  \label{einstein1}
Let $g$ and $\bar g$ be complete geodesically equivalent  metrics on a connected manifold 
$M^{n}$, $n\ge 3$. If $g$ is Einstein, then they are affine equivalent, or for a certain constants $c$, $\bar c\in \mathbb{R}$ 
the  metrics $c\cdot g$  and $\bar c\cdot \bar g$
 are Riemannian metrics of sectional  curvature $1$ (and, in particular,  the manifolds $(M^n, c\cdot g)$ and $(M^n,\bar c\cdot \bar g)$ are   finite quotients of the standard  sphere with the standard metric).     
\end{Th}

For dimension $\ge 5$, the assumption that the metrics are complete is important: if one of them is not complete, one can construct  counterexamples (essentially due to \cite{formella,mikes_einstein}). For dimensions 3 and 4, (a natural modification of)  Theorem \ref{einstein1} is  true also  locally: 
\weg{for dimension 3, it is evident since 
the metric $g$ must be homogeneous, and since by \cite{sinjukov} homogeneous metrics of nonconstant sectional curvature do not admit nontrivial geodesic equivalence. For dimension 4, we
formulate   it as }

\begin{Th}  \label{dim4}
Let $g$ and $\bar g$ be 
geodesically equivalent  metrics on a connected 3- or 4-dimensional 
manifold $M$. If $g$ is Einstein, then they are affine equivalent, or have constant sectional curvature.    
\end{Th} 

Theorem \ref{dim4}
was announced in \cite{kiosak-mikes,hinterleitner}, with the extended sketch of the proof. The proof from \cite{kiosak-mikes,hinterleitner} is very complicated: they prolonged (= covariantly differentiated)
the basic equations \eqref{basic} 6 times , and  used the condition that the metric is Einstein
at every stage of the  prolongation.  

Our proof of Theorem \ref{dim4} is a relatively easy Linear Algebra (inspired by \cite{kiosak,eastwood}) combined with a certain statement which is a relatively easy generalization of  a certain result of Levi-Civita. 

\begin{Rem} 
Theorem \ref{einstein1} is also true in dimension 2 provided the scalar curvature of $g$ is constant. Without this additional assumption  Theorem \ref{einstein1} is evidently wrong, since 
every 2-dimensional metric satisfies $R_{ij}= \tfrac{R}{2} \cdot g_{ij}.$
\end{Rem}

\subsection{History and motivation } \label{history} 

\begin{center} 
\begin{figure}[h!] \label{fig}
{{\psfig{figure=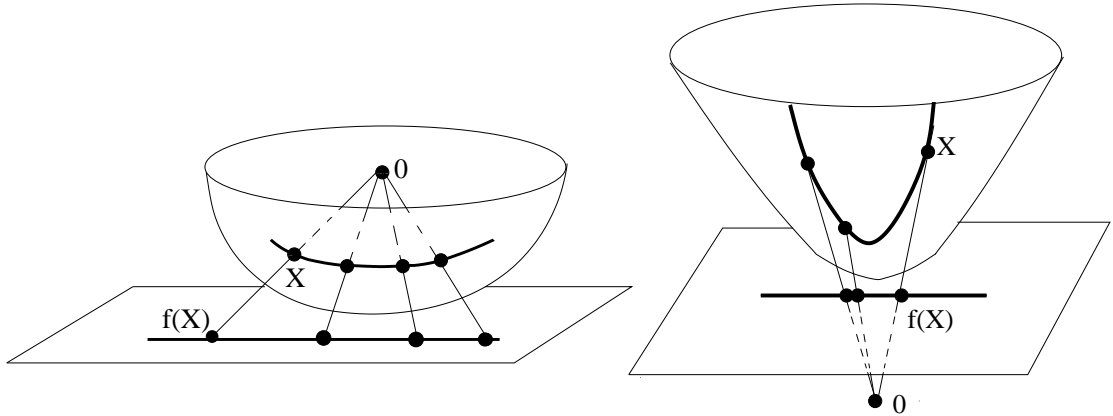}}}  
  \caption{Surfaces of constant curvature are (locally) geodesically  equivalent}\label{picture}
\end{figure}
\end{center}

The first examples of 
geodesically equivalent metrics are due to  Lagrange \cite{lagrange}. 
He observed that the radial projection $f(x,y,z)= \left(-\frac{x}{z},- \frac{y}{z}, -1\right)$  takes geodesics of the half-sphere  
$S^2:=\{(x,y,z)\in \mathbb{R}^3: \ \ x^2+y^2+z^2=1, \ z<0\}$ to the geodesics of the plane $ E^2:=\{(x,y,z)\in \mathbb{R}^3: \ \ z=-1\}$, see the left-hand side of Figure  1, since geodesics of both metrics are intersection of the 2-plane containing the point $(0,0,0)$ with the surface.  Later, Beltrami \cite{Beltrami2} generalized the example for the metrics of constant negative curvature, and for the pseudo-Riemannian metrics of constant curvature.    In the example of Lagrange, he  replaced the half sphere by  the half of 
 one of the hyperboloids $H_\pm^2:=\{(x,y,z)\in \mathbb{R}^3: \ \ x^2+y^2-z^2=\pm 1\}, $    with the restriction of the Lorenz metrics $dx^2+dy^2-dz^2$ to it. Then, the geodesics of the metric are also intersections of the 2-planes containing the point $(0,0,0)$ with the surface, and,  therefore, the stereographic projection sends it to the straight lines of the appropriate plane, see the right-hand side of Figure 1 with the (half of the) hyperboloid $H^2_-$.
  
 Though the examples of the Lagrange and Beltrami are   two-dimensional, one can easily generalize them for every dimension. 

One of the possibilities in Theorem~\ref{einstein1} is  geodesically 
equivalent metrics of constant positive  Riemannian curvature on closed manifold.  Examples of such metrics  are also  due to  Beltrami \cite{Beltrami},  we
describe their natural  multi-dimensional generalization.  Consider
the sphere $$ S^n\eqdef \{(x_1,x_2, ...,x_{n+1})\in \mathbb{R}^{n+1}: \ x_1^2+x_2^2+...+x_{n+1}^2=1\}
$$ with the metric $g $ which is the  restriction of the Euclidean metric to the sphere. Next, consider   
 the mapping $a:S^n\to S^n$  given by
$a:v\mapsto \frac{A(v)}{\|A(v)\|}$, where $A$ is an
arbitrary
non-degenerate linear  transformation of $\mathbb{R}^{n+1}$.

The mapping is clearly a diffeomorphism taking geodesics to geodesics.
Indeed, the geodesics of  $g$ are great circles
(the intersections of 2-planes that go through the origin with the
 sphere). Since $A$ is linear, it
 takes
planes to  planes. Since the normalization $w\mapsto
\frac{w}{\|w\|}$ takes punctured   planes to their intersections with the
sphere, the mapping   $a$ takes  great circles to great circles. Thus,  the pullback $a^*g$ is geodesically equivalent  to $g$.
 Evidently, if $A$ is not proportional to an
  orthogonal transformation,  $a^*g$ is not affine equivalent to $g$.

H. Weyl was probably the first who asked (in the  popular paper \cite{Weyl1}) 
whether an Einstein metric can admit   geodesically equivalent metric $\bar g$ 
 which is  nonproportional to $g$. He gave an  answer  (essentially due to \cite{Weyl2}) in the pseudo-Riemannian case 
assuming that the  metrics  $g$ and $\bar g$ have  the same light cone at every point.   Later, this question was studied by many  geometers and physicists  (a simple search in mathscinet gives about 50 papers and few books). 
 In particular, Petrov \cite{Petrov1} proved that Ricci-flat 4-dimensional Einstein metrics  of Lorenz signature can not be geodesically  equivalent, unless they are affine equivalent. It is one of the results  he obtained  the Lenin   prize (the most important  scientific award of Soviet Union) in 1972 for. 
He also explicitly asked \cite[Problem 5 on page 355]{Petrov2}   
whether the result remains true in other  dimensions. 
 
 As we will prove in Lemma \ref{einstein}, the assumption that the second metric is Einstein is not important, since it is automatically fulfilled. 
 By Theorem \ref{dim4},  the result of Petrov remains true for 4-dimensional metrics of other signatures. As we already mentioned in Section \ref{results}, the counterexamples independently  constructed by Mikes \cite{mikes_einstein}  and Formella \cite{formella} show, that the result of Petrov 
  fails in higher dimensions (so one indeed needs certain 
   additional assumptions, for example the 
  assumption that the metrics are complete   as in Theorem \ref{einstein1}, which is a standard assumption in 
 problems motivated by physics.)  
 \weg{
 The examples of \cite{mikes_einstein,formella} are warped product metrics of the form 
 \begin{equation} dx_0^2+ f(x_0)\sum_{i,j=1}^{n-1} \tilde g_{ij}(x_1,...,x_{n-1})dx_idx_j\label{warped} \end{equation}     
  where $g$ is an Einstein metric (of curvature $\tilde  R$) on a 
  $(n-1) -$dimensional manifold, and the adjusted (2-dimensional)    metric $dx_0^2+f(x_0)dy^2$ has  constant curvature   $R$.  Then, the Ricci curvature of 
 }

Resent references include  Barnes \cite{barnes},  Hall  and Lonie \cite{hall,hall3,hall4}, Hall \cite{hall1,hall2}. They studied  the existence of projective transformations of Ricci-flat metrics, which is a stronger condition than the  existence of geodesically equivalent metrics (projective transformations of $g$
allow to construct  $\bar g$ geodesically equivalent to $g$. Moreover, if $g$ is  Einstein, then $\bar g$ is  automatically Einstein as well, which essentially simplifies all  formulas). 
 
 One can find more historical details in the surveys \cite{Aminova2,eastwood1,mikes} and in the introductions to the papers \cite{topology,hyperbolic,beltrami_short,sbornik,archive,fomenko}.

{\bf Acknowledgments.} \ \
The results were obtained because 
  Gary Gibbons asked the second author to check whether certain explicitly given  Einstein metrics
  admit geodesic equivalence  (these metrics  admit integrals quadratic in velocities,
  and geodesic equivalence could lay behind  the existence of such integrals, see \cite{MT,dim2,  ERA,quantum,obata,integrable,CMH,dedicata}). 
 
 There exists  a algorithmic method  to understand whether an explicitly given metric admits an nontrivial geodesic equivalence (assuming  we can explicitly  differentiate components of the metrics, and perform algebraic operations). Unfortunately, the method  is highly computational, and applying it to the metrics suggested by Gibbons, which are given by quite complicated formulas, resulted so huge output, that we could not convince even ourself that everything is correct. 
 Therefore,  we started to look for a theory that could   simplify the calculations,  and solve the problem   in the whole generality.  
 
 We thank Gary Gibbons for his question. The second author thanks  Oxford, Cambridge,  and Loughborough Universities,     and MSRI for hospitality, and  R.~Bryant, A. Bolsinov,   
 and  M.~Eastwood for useful discussions.  Both authors were  partially supported by Deutsche Forschungsgemeinschaft
(Priority Program 1154 --- Global Differential Geometry), and by  FSU Jena.

\section{Proof of Theorem 1 \weg{\ref{einstein1}} }
 
 \subsection{Schema of the proof}

 In Section \ref{standard} we list standard facts from theory of geodesically equivalent metrics, and introduce notation we will use through the paper.
 Most of these facts can be found in the book of Sinjukov \cite{sinjukov},  but   unfortunately they are
  spread over  the text, and it in not clear under which assumption they are true (Sinjukov always assumes real-analicity, but actually needs   smoothness).   All the facts could be obtained by relatively simple tensor calculations, we will indicate how.

  The main result of Section \ref{local} are Corollaries  \ref{Tanno}, \ref{evolution}. 
In Section  \ref{proof1} we explain that the ODE along geodesics 
 given by Corollary \ref{evolution} (that controls the reparametrization that makes $g$-geodesics from $\bar g$-geodesics)
 can not have solutions such that they  satisfy the condition that both metrics are complete provided that  the Einstein  metric $g$  is pseudo-Riemannian, or  Riemannian of nonpositive scalar curvature. 
 
  Corollary \ref{Tanno}
  will be used  in Section \ref{proof2}: we will see that combining Corollary \ref{Tanno}
  with an  nontrivial result of Tanno  \cite{Tanno}  immediately gives  Theorem \ref{einstein1} under  additional assumption that the metric  is Riemannian of positive scalar curvature.

 \subsection{ Standard formulas we will use} \label{standard} 
 We work in tensor notations with the background metric $g$. That means, we sum with respect to repeating indexes, use $g$ for  raising and lowing  indexes (unless we explicitly mention), and use the Levi-Civita connection of $g$ for   covariant differentiation.

 As it was known already to Levi-Civita \cite{Levi-Civita},  two connections $\Gamma= \Gamma_{jk}^i $ and  $\bar \Gamma= \bar \Gamma_{jk}^i $  have the same unparameterized geodesics, if and only if their difference is a pure trace: there exists a $(0,1)$-tensor $\phi  $ 
such    that \begin{equation} \label{c1} 
 \bar \Gamma_{jk}^i  = \Gamma_{jk}^i + \delta_k^i\phi_{j} + \delta_j^i\phi_{k}.    
   \end{equation}

The reparameterization of the geodesics for $\Gamma$ and $ \bar   \Gamma $
connected by \eqref{c1} are done according to the following rule:  for a parametrized geodesic  $\gamma(\tau)$ of $\bar \Gamma$,  the  curve $\gamma(\tau(t)) $ is a parametrized geodesic of   $ \Gamma$,  if and only if the parameter transformation $\tau(t)$ satisfies the following ODE: 
\begin{equation} \label{umparametrisation} 
\phi_\alpha \dot \gamma^\alpha =  
\frac{1}{2} \frac{d}{dt} \left(\log\left(\left|\frac{d\tau}{dt}\right|\right)\right). 
\end{equation} 
(We denote by $\dot \gamma$ the velocity vector of  $\gamma$ with respect to the parameter $t$,  and assume summation with respect to repeating index $\alpha$.) 

If $\Gamma$  and   $\bar \Gamma$  related by  \eqref{c1}  are Levi-Cevita connections of  metrics $g$ and $\bar g$, then one can find explicitly (following Levi-Civita \cite{Levi-Civita}) a function $\phi$ on the manifold such that its differential $\phi_{,i}$  coincides with the covector $\phi_i$: indeed, contracting \eqref{c1}  with respect to  $i$ and $j$, we obtain 
 $\bar \Gamma_{\alpha i}^\alpha   = \Gamma_{\alpha i}^\alpha  + (n+1) \phi_{i}$. 
  From the  other side, for the Levi-Civita 
 connection  $\Gamma$ of a metric $g$  we have  $  \Gamma_{\alpha k}^\alpha  = \tfrac{1}{2} \frac{\partial \log(|det(g)|)}{\partial x_k} $.  Thus, 
 \begin{equation} \label{c1,5}  \phi_{i}= \frac{1}{2(n+1)}  \frac{\partial }{\partial x_i}  \log\left(\left|\frac{\det(\bar g)}{\det( g)}\right|\right)= \phi_{,i} \end{equation} 
  for the function $\phi:M\to \mathbb{R}$ given by 
  \begin{equation} \label{phi}  \phi:= \frac{1}{2(n+1)} \log\left(\left|\frac{\det(\bar g)}{\det( g)}\right|\right). \end{equation}  In particular, the derivative of $\phi_i$ is  symmetric, i.e., $\phi_{i,j}= \phi_{j,i}$.

The formula  \eqref{c1}   implies  that two metrics $g$ and $\bar g$ are geodesically equivalent if and only if   for a certain $\phi_{i}$ (which is, as we explained above, the differential of $\phi$ given by  \eqref{phi}) we have 
\begin{equation}\label{LC}
    \bar g_{ij, k} -  2 \bar g_{ij} \phi_{k}-  \bar g_{ik}\phi_{j} -   \bar g_{jk}\phi_{i}= 0, 
\end{equation} 
where ``comma" denotes the covariant derivative with respect to the connection $\Gamma$. 
Indeed, the left-hand side of this equation is the covariant derivative with respect to  $\bar \Gamma$, and vanishes if and only if  $\bar \Gamma$ is the Levi-Civita connection for $\bar g$.

The equations \eqref{LC} can be linearized by a clever  substitution:  consider 
   $a_{ij}$ and $\lambda_i$ given by 
\begin{eqnarray} \label{a}
a_{ij} &=   &e^{2\phi} \bar g^{\alpha \beta} g_{\alpha i} g_{\beta j},\\  \label{lambda} 
\lambda_{i} & = &  -e^{2\phi}\phi_\alpha \bar g^{\alpha \beta} g_{\beta i}, \end{eqnarray}
where $\bar {g}^{\alpha \beta}$ is  the tensor dual to $\bar g_{\alpha \beta}$:  \ $\bar {g}^{\alpha i}\bar g_{\alpha j}= \delta_j^i$. 
It is an easy exercise to show that the following linear  equations on  the symmetric $(0,2)-$tensor $a_{ij}$ and $(0,1)-$tensor $\lambda_i$ are    equivalent to \eqref{LC}. 
 \begin{equation} \label{basic} 
 a_{ij,k}= \lambda_i g_{jk} + \lambda_j  g_{ik}. 
 \end{equation}

 \begin{Rem} 
 For dimension 2, the substitution  (\ref{a},\ref{lambda})  was already known to R.  Liouville \cite{liouville} and Dini \cite{Dini}, see \cite[Section 2.4]{bryant} for details and a conceptual explanation.  For arbitrary dimension,  the substitution (\ref{a},\ref{lambda}) and the equation \eqref{basic} are due to Sinjukov \cite{sinjukov}.    Geometry staying behind  is explained in \cite{eastwood}.
 \end{Rem}

Note that it is possible to find a  function $\lambda$ such that its differential is precisely 
the the $(0,1)-$tensor $\lambda_i$: indeed, multiplying \eqref{basic}  by $g^{ij}$ and summing with respect to repeating indexes $i,j$ we obtain $(g^{ij}a_{ij})_{,k} = 2  \lambda_k$. Thus,
$\lambda_i$ is the differential of the function 
\begin{equation}\label{lam} 
\lambda:= \tfrac{1}{2} g^{\alpha \beta }a_{\alpha \beta}.   
\end{equation}  
In particular, the covariant derivative of $\lambda_i$ is symmetric:  
 $\lambda_{i,j} = \lambda_{j,i}$.

 Integrability conditions for the equation \eqref{basic}  (we substitute the derivatives of $a_{ij}$ given by \eqref{basic} in the formula  $a_{ij,lk}- a_{ij,kl}= a_{i \alpha }R^{\alpha}_{jkl} +  a_{\alpha j}R^{\alpha}_{ikl}$, which is true for every $(0,2)-$tensor  $a_{ij}$)   were  first obtained by Solodovnikov \cite{s1} and are   
 
\begin{equation}  a_{i \alpha }R^{\alpha}_{jkl} +  a_{\alpha j}R^{\alpha}_{ikl} =\lambda_{ l,i} g_{jk}+\lambda_{ l,j} g_{ik}-\lambda_{ k,i} g_{jl}-\lambda_{ k,j} g_{il}. \label{int1} \end{equation}

For further use let us recall the fact which can also be obtained  by simple calculations:  the 
Ricci-tensors of connections related by \eqref{c1} are  connected by the formula 
\begin{equation} \label{ric} \bar R_{ij} = R_{ij}- (n-1)(\phi_{i,j}- \phi_{i}\phi_{j}),  
\end{equation} 
 where $R_{ij}$ is the Ricci-tensor  of $\Gamma  $ and $\bar R_{ij}$ is the Ricci-tensor  of $\bar \Gamma $.

\subsection{Local results} \label{local} 
Within the whole paper we work on a smooth manifold of dimension $n\ge 3$.

\begin{Lemma}[Folklore] \label{first} 
Let $a_{ij}$ be a solution of \eqref{basic}
for the  metric $g$. Then, it commutes with the Ricci-tensor: \begin{equation} \label{en0}  
a^{\alpha}_{i}R_{\alpha j} = a^{\alpha}_{j}R_{i \alpha }. \end{equation} 
\end{Lemma} 

{\bf Proof. }  Consider the equations \eqref{int1}.   We ``cycling" the equation with respect to $i,k,l$: we sum it with itself after renaming the indexes according to $(i\mapsto k \mapsto l \mapsto i)$ and 
with  itself after renaming the indexes according to $(i\mapsto l \mapsto k \mapsto i)$. 
The first term at the left-hand side of the equation will disappear because of the Bianchi equality $R^{\alpha}_{ikl} +  R^{\alpha}_{kli} + R^{\alpha}_{lik}=0 $, the right-hand side vanishes completely,  
 and we obtain 

\begin{equation} 
a_{\alpha i} R^\alpha_{jkl} + a_{\alpha k} R^\alpha_{jli}+ a_{\alpha l} R^\alpha_{jik}=0.
\end{equation} 

Multiplying with $g^{jk}$, using the symmetries of the curvature tensor,   and summing over the repeating indexes we obtain $a_{\alpha i} R^\alpha_l - a_{\alpha l} R^{\alpha }_i=0 $ implying the claim, \qed 

\begin{Lemma} \label{harmonic}    Suppose the  curvature tensor of the metric $g $ satisfies 
$$
R^{\alpha}_{ijk, \alpha}=0.$$
Then, for every solution $a_{ij}$  of \eqref{basic} 
such that $\lambda_i\ne 0$ at a point $p\in M^n$, in a sufficiently  
small neighborhood $U(p)$ of $p$ we have  
 \begin{equation} \label{goal} 
  \lambda_{k,j} = \C{1} g_{kj} + \C{2} R_{kj}+ \C{3} a_{kj}+ \C{4} a^{\alpha}_jR_{\alpha k},   
\end{equation} 
where the coefficients $\C{1}$, $\C{2}$, $\C{3}$, $\C{4}$  are given by the formulas 
$$\C{1}=    \frac{-\lambda_\alpha a^{\alpha}_\beta  \xi^\beta R + 2 \lambda \lambda ^{\ \ \ \beta}_{ \alpha,\beta} \xi^\alpha + a^\alpha_\beta R^\beta_\alpha - 4\lambda_{\alpha,}^{\ \  \alpha}}{4n}; \ \ 
 \C{2}=   \tfrac{1}{4}  \lambda_\alpha a^{\alpha}_\beta \xi^\beta ;   \  \ 
 \C{3}= - \tfrac{1}{4}\lambda^{\ \ \ \beta}_{\alpha, \beta} \xi^\alpha; \C{4} =  - \tfrac{1}{4}a^\alpha_\beta   R_{\alpha}^\beta  
 , 
$$where $\xi$ is an arbitrary vector field such that $\lambda_i\xi^i=1$. 
\end{Lemma}

\begin{Rem} \label{nov} 
The assumptions  of the lemma are automatically fulfilled for Einstein spaces. Indeed, the second  Bianchi  identity for the curvature tensor is 
$$
R^{h}_{ijk,l} +  R^{h}_{ikl,j}+ R^{h}_{ilj,k}=0.
$$
Contracting with respect to $h$ and  $l$, we obtain 
$$
R^{\alpha}_{ijk,\alpha} +  \underbrace{R^{\alpha }_{ik\alpha,j}}_{-R_{ik,j}}+ \underbrace{R^{\alpha}_{i\alpha j,k}}_{R_{ij, k} }=0.
$$
If the metric is Einstein, then the second and the third components of the equation vanishes, and we obtain 
$R^{\alpha}_{ijk,\alpha}=0$. Moreover, we see that  actually the condition $R_{ik,j} - R_{ij, k}=0$ is a necessary and  sufficient condition  for      $R^{\alpha}_{ijk,\alpha}=0$. \end{Rem}   \begin{Rem} 
The tensor $R_{ik,j} - R_{ij, k}$ is called projective Yano tensor, and plays important role in the theory of geodesically equivalent metrics; in particular, it is projectively invariant in dimension 2 \cite{liouville,bryant}, and is  an  essential   part of the so-called tractor approach for the investigation of geodesically equivalent metrics \cite{eastwood}. 
\end{Rem}

{\bf Proof of Lemma \ref{harmonic}.} 
Consider the solution $a_{ij}$  of the equation 
\eqref{basic}.  Let us take the covariant derivative of the equations \eqref{int1} (the index of differentiation is ``$m$"), and replace the covariant derivative of $a$ by formula \eqref{basic}.  We obtain 

\begin{equation} \label{en1} 
\begin{array}{l} \  \  \ \lambda_\alpha R^{\alpha} _{jkl} g_{im} + \lambda_i R_{mjkl} + a_{\alpha i}  R^{\alpha}_{jkl, m} +
\lambda_\alpha R^{\alpha} _{ikl} g_{jm} + \lambda_j R_{mikl} + a_{\alpha j}  R^{\alpha}_{ikl, m}\\ =\lambda_{ l,im} g_{jk}+\lambda_{ l,jm} g_{ik}-\lambda_{ k,im} g_{jl}-\lambda_{ k,jm} g_{il}.\end{array}
\end{equation} 
  
We multiply with $g^{lm}  $,  sum  with respect to repeating indexes $l,m$, and use $
R^{\alpha}_{ijk, \alpha}=0.$  We obtain: 

\begin{equation} \label{en2} 
 \lambda_\alpha R^{\alpha} _{ikj}  + \lambda_\alpha R^{\alpha} _{jki} - \lambda_i R_{jk}- \lambda_j R_{ik}
 =\lambda^{\ \ \ \alpha}_{ i, \alpha} g_{jk}+\lambda^{\ \ \ \alpha}_{ j, \alpha} g_{ik}
 -\lambda_{ k,ij} -\lambda_{ k,ji}.
\end{equation} 
We now alternate the equation \eqref{en2}  with respect to $k,j$ to obtain
 
\begin{equation} \label{en3} 
\begin{array}{l} \  \  \ 4 \lambda_\alpha R^{\alpha} _{ikj}   
 =\lambda^{\ \  \ \alpha}_{ j, \alpha} g_{ik}
 -\lambda^{\ \ \  \alpha}_{ k, \alpha} g_{ij}- \lambda_k R_{ij} + \lambda_j R_{ik}.\end{array}
\end{equation}

Let us now rename the indexes $i\mapsto k\mapsto j \mapsto \alpha$ in   \eqref{en3},  multiply the result by $a^\alpha_i$, use the symmetries of the curvature tensor  and sum over the  repeating index $\alpha$. We obtain

\begin{equation} \label{en4} 
\begin{array}{lcl} \  \  \ 4  a^\alpha_i  R_{\alpha jk\beta} \lambda_{}^{\beta} 
&=&   4  a^\alpha_i  R_{k j \alpha }^\beta \lambda_{\beta} \\ 
&=&  a^\alpha_i \left(  \lambda^{\  \  \  \beta}_{  \alpha, \beta} g_{kj}
 -\lambda^{\  \  \  \alpha}_{ k, \alpha} g_{ij}  -\lambda^{\ \ \ \beta}_{  j, \beta} g_{k\alpha}
 - \lambda_j R_{\alpha k} + \lambda_\alpha  R_{ j k}\right) \\ 
 &=&  a^\alpha_i   \lambda^{\ \ \ \beta}_{  \alpha, \beta} g_{kj} 
 - \lambda^{\ \ \ \beta}_{j,\beta} a_{ki }
- \lambda_j a^\alpha_i R_{\alpha k} + \lambda_\alpha  a^{\alpha}_i R_{kj} .
  \end{array}
\end{equation} 

Now we multiply the equation \eqref{int1} by  $\lambda^l$ and sum over the repeating index $l$. We see that the first component of the result is precisely the left-hand side of \eqref{en4}; we replace it by 
the right-hand side  of \eqref{en4}. We obtain

\begin{equation} \label{en51} \begin{array}{lll} 
0& =& \left( a^\alpha_i   \lambda^{\ \ \ \beta}_{   \alpha, \beta}- 4 \lambda^{ \alpha}\lambda_{\alpha, i} \right)  g_{kj} - \lambda^{\ \ \ \beta}_{j,\beta} a_{ki }
+ \lambda_j \left(  - a^\alpha_i R_{\alpha k} + 4 \lambda_{k,i}\right) + \lambda_\alpha  a^{\alpha}_i R_{kj} \\
&+& \left( a^\alpha_j   \lambda^{\  \ \  \beta}_{    \alpha, \beta}- 
4\lambda^{\alpha}\lambda_{ \alpha, j} \right)  g_{ki} - \lambda^{\  \  \ \beta}_{i,\beta} a_{kj }
+ \lambda_i \left( - a^\alpha_j R_{\alpha k} + 4 \lambda_{k,j}\right) + \lambda_\alpha  a^{\alpha}_j R_{ki}
\end{array} \end{equation}

We now  alternate  \eqref{en4} 
with respect to $k,j$,  rename $k \longleftrightarrow i$, and add  the result to  \eqref{en51}. After dividing by 2 for cosmetic reasons,  and using that by Lemma   \ref{first} the tensor   $a^\alpha_i R_{\alpha k}$  is symmetric  with respect to $i,k$,   we obtain

\begin{equation} \label{en5}
\left( a^\alpha_i   \lambda^{\ \ \ \beta}_{  \alpha, \beta} - 4 \lambda^{\alpha}\lambda_{ \alpha, i} \right)  g_{kj} 
 + \lambda_\alpha  a^{\alpha}_i R_{kj} 
 - \lambda^{ \  \ \ \beta}_{i,\beta} a_{kj }
+ \lambda_i \left( - a^\alpha_j R_{\alpha k} +4  \lambda_{k,j}\right) =0.
\end{equation} 

We multiply \eqref{en5} by $g^{kj}$  and sum over the  repeating indexes $k,j$. We obtain (after dividing by $n$) 
\begin{equation} \label{en6}
\left( a^\alpha_i   \lambda^{ \ \ \  \beta}_{  \alpha, \beta} - 4 \lambda^{\alpha}\lambda_{ \alpha,  i} \right) = 
- \tfrac{R}{n}  \lambda_\alpha  a^{\alpha}_i 
+  \frac{2\lambda}{n}  \lambda^{ \  \ \ \beta}_{ i,\beta} 
- \lambda_i \frac{\left(-  a^\alpha_\beta  R_{\alpha}^{ \beta} +4  \lambda_{\alpha,  }^{\ \  \alpha}\right)}{n} =0,
\end{equation} 
where $R:=R_{\alpha \beta}g^{\alpha \beta}$ is the scalar curvature of $g$.  Substituting the expression for $\left( a^\alpha_i   \lambda^{\ \ \ \beta}_{  \alpha, \beta} - 4 \lambda^{\alpha}\lambda_{ \alpha, i} \right)$ from \eqref{en6}  in \eqref{en5}, we obtain

\begin{equation} \label{en7}
\begin{array}{ccl} 0&=& \lambda_\alpha a^{\alpha}_{  i}\left( R_{kj} - \tfrac{R}{n}g_{kj} \right) + \lambda_{i,\beta}^{\ \ \ \beta} \left( \frac{2\lambda}{n} g_{kj} - a_{kj} \right) \\ & -& \lambda_i\left(  \frac{-a^\alpha_\beta  R_{\alpha}^{ \beta} +4  \lambda_{\alpha, }^{\ \  \alpha}}{n}g_{kj} + a^{\alpha}_j R_{\alpha k} - 4 \lambda_{k,j}\right)\end{array} 
\end{equation} 

Since $\lambda_{i}\ne 0$  at a point $p$,  then  $\lambda_i\xi^i=1$ for a  certain vector field $\xi$ in  a sufficiently small neighborhood $U(p)$.  Contracting the equation  \eqref{en7} with this $\xi^i$,  we obtain  

\begin{equation} \label{en8}\begin{array}{ccl} 0&=& 
\lambda_\alpha a^{\alpha}_{  i} \xi^i \left( R_{kj}  -\tfrac{R}{n}g_{kj} \right) + \xi^i\lambda_{i,\beta}^{\ \ \beta} \left( \frac{2\lambda}{n} g_{kj} - a_{kj} \right) \\ &-& \lambda_i\left(  \frac{-a^\alpha_\beta  R_{\alpha}^{\beta} +4  \lambda_{\alpha, }^{\ \ \ \alpha}}{n}g_{kj} + a^{\alpha}_j R_{\alpha k} - 4 \lambda_{k,j}\right)\end{array} 
\end{equation} 

We see that $\lambda_{j, k}$ is a linear combination of $a^{\alpha}_j R_{\alpha k}$, $g_{jk}$, $R_{jk}$ and $a_{kj}$  as we want. The coefficients in the linear combination are as in the formula  below , \qed  
$$ 
 4 \lambda_{k,j}=   a_{\alpha k}   R^{\alpha}_j + 
 \frac{-\lambda_\alpha a^{\alpha}_\beta  \xi^\beta R + 2 \lambda \lambda ^{\ \ \ \beta}_{ \alpha,\beta} \xi^\alpha + a^\alpha_\beta R^\beta_\alpha - 4\lambda_\alpha^\alpha}{n} g_{jk} +   \lambda_\alpha a^{\alpha}_\beta \xi^\beta  R_{jk}  - \lambda^{\ \ \ \beta}_{\alpha, \beta} \xi^\alpha  a_{kj} . 
$$

 \begin{Cor} 
 Assume $g$ is an Einstein metric. Let $a_{ij}$ be a solution of \eqref{basic}.  Assume $\lambda_{i}\ne 0$ at a point $p$. Then, in a  sufficiently small   neighborhood of $p$,   
  $\lambda_{i,j}$ is a linear combination of $g_{ij}$ and  $a_{ij}$: 
\begin{equation}\label{vb} 
 \lambda_{i,j}= \mu g_{ij}+ K a_{ij},  
 \end{equation}
 where  the coefficients  $K:=-\tfrac{R}{n(n-1)}$ and  $\mu:= \frac{\lambda_{\alpha,}^{\ \alpha}- 2 K \lambda}{n}$.
 \end{Cor}
 
 {\bf Proof.}  By assumption, in a small neighborhood of $p$ we have  $\lambda_i\ne 0 $; 
 this implies that $a_{ij}$ is not proportional to $g_{ij}$, because  by the result of Weyl \cite{Weyl2}
 if two metrics are geodesically and conformally equivalent, then  they are proportional (with a constant coefficient of proportionality). 
 
  As we explained in Remark  \ref{nov},  
 the assumptions of the lemma are fullfilled if the metric is Einstein. Moreover, if the   metric is Einstein,  then the second term of the right-hand side of \eqref{goal} is proportional to $g$, and the last term is proportional to $a$ implying that $\lambda_{i,j}$ is a linear combination of $g_{ij}$ and $a_{ij}$. We need to calculate the coefficients of the linear combination.  
 
 Substituting  the condition that the metric is Einstein in \eqref{en3} we obtain   
 \begin{equation} \label{cor1} 
\begin{array}{l} \lambda_\alpha R^{\alpha}_{ikj}   
 =\tau_j g_{ik} - \tau_k g_{ij}
 ,\end{array}
\end{equation} 
 where 
 \begin{equation} \label{cor2} \tau_i:= \tfrac{1}{4} \left(\lambda_{i,\alpha}^{ \  \   \ \alpha}+ \tfrac{R}{n}\lambda_i\right).\end{equation}
 Contracting the equation \eqref{cor2} with  $g^{ij}$ we obtain  $(n-1)\tau_j=  -\tfrac{R}{n} \lambda_j$ implying  
 \begin{equation} \label{cor3} \tau_j=  -\tfrac{R}{n(n-1)} \lambda_j.\end{equation}
 Now, since the metric is Einstein, the first bracket  in the 
 sum  \eqref{en7} is zero, and the term $a^\alpha_\beta  R_{\alpha}^{ \beta}$ equals  $\tfrac{R}{n}\delta_{\alpha}^{ \beta}a^\alpha_\beta =2\tfrac{R}{n} \lambda$,   so the formula \eqref{en7} reads 
 
\begin{equation} \label{cor4} \lambda_{i,\beta}^{\ \ \ \beta} \left( \frac{2\lambda}{n} g_{kj} - a_{kj} \right) - \lambda_i\left(  \frac{-  2 \lambda \tfrac{R}{n}  +4  \lambda_{\alpha, }^{\  \alpha}}{n}g_{kj} + \tfrac{R}{n}a_{kj}   - 4 \lambda_{k,j}\right)=0
\end{equation} 
 
Combining \eqref{cor2} and \eqref{cor3},  we obtain \begin{equation} \label{new1}\lambda_{i,\beta}^{\ \ \ \beta}= \left(4 k- \tfrac{R}{n} \right)\lambda_i.\end{equation}  Substituting this in \eqref{cor4}, we obtain 
\begin{equation} \label{cor5} \left(4 k- \tfrac{R}{n} \right) \left( \frac{2\lambda}{n} g_{kj} - a_{kj} \right)\lambda_i -\left(  \frac{ - 2 \lambda \tfrac{R}{n}  +4  \lambda_{\alpha, }^{ \ \alpha}}{n}g_{kj} + \tfrac{R}{n}a_{kj}   - 4 \lambda_{k,j}\right)\lambda_i=0.
\end{equation} 
Since by assumption $\lambda_i\ne 0$, we obtain \eqref{vb}, 
 \qed

 \begin{Rem} 
 Assume $g$ is an Einstein metric. Let $a_{ij}$ be a solution of \eqref{basic}.  Then, 
\begin{equation}\label{yano} 
 \lambda_{\alpha} Y^\alpha_{ijk}= 0,  
 \end{equation}
 where $  Y^h_{ijk}:= R^h_{ijk} - \tfrac{R}{n(n-1)}\left(\delta^h_j g_{ik}-\delta^h_kg_{ij} \right)$ is the so-called 
   concircular curvature of $g$ introduced by Yano \cite{Yano}.  
 \end{Rem} 
 
 {\bf Proof. } Substituting \eqref{cor3} in \eqref{cor1},  we obtain the claim, \qed
 
  \begin{Cor} \label{vbcor} 
 Assume $g$ is an Einstein metric. Let $a_{ij}$ be a solution of \eqref{basic}.  Consider $K:=-\tfrac{R}{n(n-1)}$ and the function $\mu:= \frac{\lambda_{\alpha,}^{\ \ \alpha}- 2 K \lambda}{n}$. Then, the function $\mu$ satisfies  the equation 
\begin{equation}\label{vb2} 
 \mu_{,i}= 2  K \lambda_{i}.  
 \end{equation}
 \end{Cor}
 \begin{Rem}
In particular, under the assumptions of   Corollary \ref{vbcor}, for a certain 
  $\const\in \mathbb{R}$,   the function 
  $\lambda + \const$  is an eigenfunction of the laplacian of $g$. 
 \end{Rem} 
 {\bf Proof of Corollary \ref{vbcor}.} If $\lambda$ is constant  in a neighborhood of a point, the equation \eqref{vb2} is automatically fulfilled. Below  we will  assume that $\lambda$ is not constant. Differentiating the definition of $\mu$ and multiplying by $n$ for cosmetic reasons, we obtain 
 \begin{equation} \label{new0} 
 n\mu_{,i}=2 \lambda_{\alpha ,\ i}^{\ \  \alpha} -2 K \lambda_{i}.   
\end{equation}
 By definition of curvature we have $\lambda_{i,jk}-\lambda_{i,kj}= \lambda_\alpha R^{\alpha}_{ijk}$.
 Contracting this  with $g^{ij}$, and using $R_{ij}= \tfrac{R}{n}g_{ij}$, we obtain 
 $$
 \lambda_{\alpha,  \  \ k}^{\ \   \alpha} - \lambda_{k,\alpha }^{ \ \ \ \alpha}  = -\tfrac{R}{n}\lambda_{k}. 
 $$
   The formula \eqref{new1} gives us $\lambda_{k,\alpha }^{\ \ \ \alpha} $, whose substitution gives      $$
  \lambda_{\alpha, \ k}^{\ \ \alpha}= \left(-2 \tfrac{R}{n} + 4 K \right)\lambda_{k}.
  $$ Substituting this in \eqref{new0}, we obtain 
  $\mu_{,i} =-\tfrac{2 R}{n(n-1)}\lambda_{i}= 2 K\lambda_{i},$ \qed  
  
\begin{Cor} \label{Tanno} 
Let $g$ and $\bar g$ be geodesically equivalent metrics, assume $g$ is an Einstein metric. Then,  the function $\lambda$ given by \eqref{lam} satisfies  
   \begin{equation} \label{tanno} 
\lambda_{,ijk} - K\cdot \left(2 \lambda_{,k} g_{ij} + \lambda_{,j} g_{ik} + \lambda_{,i}g_{jk}\right)=0, \end{equation} 
 where  $K:=-\tfrac{R}{n(n-1)}$. 
   \end{Cor} 
   
   {\bf Proof. }  If $\lambda$ is constant in a neighborhood of $p$, the equation is automatically fulfilled.    Then, it is sufficient to prove Corollary \ref{Tanno}  at points $p$ such  that $\lambda_{i}(p)\ne 0$. 
   
Covariantly    differentiating \eqref{vb}, we obtain $\lambda_{i,jk}= \mu_{,k}g_{ij} +  K a_{ij,k}$. Substituting $\mu_{,k}$ by \eqref{vb2}, and $a_{ij,k}$ by \eqref{basic}, we obtain the claim, \qed

 \begin{Lemma} \label{einstein}
 Let $g$ and $\bar g$ be  geodesically equivalent. Assume $g$ is Einstein, and assume that 
  $\lambda_{i}\ne 0$ at a point $p$.   
 
    Then, the restriction of $\bar g$ to a sufficiently small neighborhood $U(p)$ is Einstein as well. Moreover, the following formula holds (at every point of  $U(p)$).
 \begin{equation}\label{f1} 
 \phi_{i,j} - \phi_{i} \phi_{j} = \tfrac{R}{n(n-1)} g_{ij}-\tfrac{\bar R}{n(n-1)} \bar g_{ij}  ,
 \end{equation} 
 where $\bar R$ is the scalar curvature of the metric $\bar g$.
 \end{Lemma}
 
 {\bf Proof. }   
 We covariantly  differentiate \eqref{lambda} (the index of differentiation is ``j"); then we substitute the expression \eqref{LC} for $\bar g_{ij,k}$   to obtain 
 \begin{equation} \label{f2} \begin{array}{ccl}
 \lambda_{i,j} &=& -2 e^{2\phi}\phi_{j} \phi_\alpha \bar g^{\alpha \beta} g_{\beta i}-e^{2\phi}\phi_{\alpha,j} \bar g^{\alpha \beta} g_{\beta i}+e^{2\phi}\phi_\alpha  \bar g^{\alpha \gamma} \bar g_{\gamma l,j} \bar g^{l \beta} g_{\beta i} \\ &=&   -e^{2\phi}\phi_{\alpha,j} \bar g^{\alpha \beta} g_{\beta i}+e^{2\phi}\phi_\alpha \phi_\gamma \bar g^{\alpha \gamma} \bar  g_{ i j }+  e^{2\phi}\phi_{j} \phi_l \bar g^{l\beta}g_{\beta i}      \end{array}    ,  
 \end{equation} 
 where $\bar g ^{\alpha \beta}$ is the tensor dual to $\bar g_{\alpha \beta}$.  
 We now  substitute  $\lambda_{i,j}$ from \eqref{vb},  use that $a_{ij}$ is given by \eqref{a}, and divide by $e^{2\phi}$ for cosmetic reasons   to  obtain 
 \begin{equation} \label{f3} 
 e^{-2\phi} \mu g_{ij} + K \bar g^{\alpha \beta} g_{\alpha j}g_{\beta _i} = -\phi_{\alpha, j} \bar g^{\alpha \beta} g_{\beta i}+\phi_\alpha \phi_{\gamma} \bar g^{\alpha \gamma} \bar  g_{ i j }+ \phi_{j} \phi_{l} \bar g^{l\beta}g_{\beta i}.  
 \end{equation}  
 Multiplying with $g^{i\xi} \bar g_{\xi k}$,  we obtain 
 \begin{equation} \label{f4} 
 \phi_{k,j}-\phi_{k}\phi_{j} =(\phi_\alpha \phi_\beta \bar g^{\alpha \beta} - e^{-2 \phi } \mu ) \bar g_{kj} - K g_{kj}.  
  \end{equation} 
  Let us now show that the coefficient 
  ${\bar K}{} := -\tfrac{\phi_\alpha \phi_{\beta} \bar g^{\alpha \beta} - e^{-2 \phi } \mu}{n-1}$ is constant.  Substituting 
  \eqref{f4} in \eqref{ric}, and using $R_{ij}= \tfrac{R}{n} g_{ij}$, we obtain 
  $$
  \bar R_{ij} = \tfrac{R}{n} g_{ij} - \tfrac{R}{n} g_{ij} - (\phi_{\alpha} \phi_{\beta} \bar g^{\alpha \beta} - e^{-2 \phi } \mu ) \bar g_{ij}. 
  $$ 
  We see that $\bar R_{ij} $ 
  is proportional to $\bar g_{ij}$. Then,  $\bar g$ is an Einstein metric; in particular, $\bar K:=-\tfrac{ \phi_{\alpha} \phi_{\beta} \bar g^{\alpha \beta} - e^{-2 \phi } \mu }{n-1}$ is a constant equal to $ -\tfrac{\bar R}{n(n-1)}$, and \eqref{f4} gives us  the formula 
 \begin{equation} 
 \bar K \bar g_{ij}= K g_{ij}+ \phi_{i,j} - \phi_{i}\phi_{j}, 
\label{phi1} 
 \end{equation} which 
 is evidently equivalent to \eqref{f1}, \qed

\begin{Cor} \label{evolution}
Let $g$ and $\bar g$ be geodesically equivalent metrics, assume $g$ is an Einstein metric. Consider a (parametrized)  geodesic $\gamma$ of the metric $g$, 
 and denote by $\dot \phi$, $\ddot \phi$ and $\dddot \phi$ the first, second and third derivatives of the function $\phi$  given by \eqref{phi} along the geodesic. Then, along the geodesic,   
 the following ordinary differential  equation holds: 
   \begin{equation}  
 \begin{array}{lcl}\dddot\phi&=&  4  K g(\dot\gamma, \dot\gamma) \dot\phi  + 6 \dot\phi\ddot\phi- 4(\dot \phi)^3\end{array}, 
  \label{phi0}
 \end{equation} where  $g(\dot\gamma, \dot\gamma):=  g_{ij} \dot\gamma^i\dot\gamma^j$. 
   \end{Cor} 
  {\bf Proof. }  If $\phi\equiv 0$ in a neighborhood $U$,  the equation is automatically fulfilled.  
 Then, it is sufficient to prove Corollary \ref{evolution} assuming $\phi_{i}$ is not constant. 
 
  The formula \eqref{f1} is evidently equivalent to  \eqref{phi1}, which is evidently equivalent to 
 \begin{equation}
 \phi_{i,j}=\bar K \bar g_{ij}- K g_{ij}+   \phi_{i}\phi_{j}. 
 \label{phi2} 
 \end{equation}
Taking covariant derivative of  \eqref{phi2},  we obtain 
  \begin{equation}
 \phi_{i,jk}=\bar K \bar g_{ij,k} +    2 \phi_{i,k}\phi_{j}+ 2 \phi_{j,k}\phi_{i}. 
  \label{phi3} 
 \end{equation}
 Substituting the expression for $\bar g_{ij,k} $ from \eqref{LC}, and  substituting $\bar K \bar g_{ij}$ given by  
\eqref{phi1}, we obtain  
 \begin{equation}
 \begin{array}{lcl}\phi_{i,jk}&=&\bar K  ( 2 \bar g_{ij} \phi_{k}+   \bar g_{ik}\phi_{j} +    \bar g_{jk}\phi_{i})+    2 \phi_{i,k}\phi_{j}+ 2 \phi_{j,k}\phi_{i} \\&=& K(  2  g_{ij} \phi_{k}+   g_{ik}\phi_{j} +   g_{jk}\phi_{i} ) + 2 (\phi_{k}\phi_{i,j} + \phi_{i}\phi_{j,k}+\phi_{j}\phi_{k,i})- 4\phi_{i}\phi_{j}\phi_{k} \end{array}
  \label{phi4}
 \end{equation}
 Contracting with $\dot\gamma^i \dot\gamma^j \dot\gamma^k$ and using that $\phi_i$ is the differential of the function \eqref{phi}  we obtain the desired ODE \eqref{phi0}, \qed 

 \begin{Cor} Let $\bar g$ (on a connected $M^{n\ge 3})$ \label{100}
 be geodesically equivalent to an Einstein metric $g$, but is not affine equivalent ot $g$. Then,  
 the restrictions of $g$ and $\bar g$ to any neighborhood are also not affine equivalent. 
 \end{Cor}  
 
 \begin{Rem} 
 The assumption that $g$ is  Einstein is important: Levi-Civita's description of geodesically equivalent metrics \cite{Levi-Civita}   immediately gives  
 counterexamples.  
 \end{Rem} 
 
 {\bf Proof of Corollary \ref{100}.} Suppose  $\phi_i\ne 0$ at a point $p$. Consider a geodesic $\gamma(t)$  such that $\gamma(0)=p$, $\dot\gamma^\alpha(0) \phi_\alpha \ne 0 $. By Lemma \ref{evolution}, the function $\phi_i \dot\gamma^i$ satisfies equation \eqref{phi0} along the geodesic. Since $\dot \phi(0)\ne 0$, then $\dot \phi(t) \ne 0 $ for almost every $t$.
 Then,   $\phi_i \ne 0$ at almost every point of geodesic. Since every point can be reached from the point $p$ by a sequence of geodesics, $\phi_i\ne 0$   at almost every point, \qed 
 
  \subsection{ Proof of Theorem  \ref{einstein1} for Riemannian metrics of nonpositive  scalar  curvature, and for pseudo-Riemannian metrics}
  
  Assume the metric $g$ on a connected $M^{n\ge 3}$ 
   is Einstein and is either Riemannian (i.e., positive defined) with nonpositive  scalar curvature,
    or essentially  pseudo-Riemannian (i.e., there exists light-like vectors). 
    Let $\bar g$ be geodesically equivalent to $g$. Assume  both metrics are complete. 
 Our goal is to    show that $\phi$ given by \eqref{phi} is constant, because  in view of \eqref{c1}  this  implies that the metrics are affine equivalent.
 
 Consider a parameterized  geodesic $\gamma(t)$ of $g$.  If the metric $g $  is pseudo-Riemannian, we additionally  assume that $\gamma$ is a light-like geodesic i.e., $\dot \gamma^i  \dot \gamma^j g_{ij}=0$. 
  Since the metrics are geodesically equivalent, for a certain function $\tau:\mathbb{R}\to \mathbb{ R}$  the 
   curve $\gamma(\tau)$ is a geodesic of $\bar g$.     
 Since the metrics are complete, 
 the reparameterization  $\tau(t)$ is a diffeomorphism $\tau:\mathbb{R}\to \mathbb{R}$. Without loss of generality we can think that $\dot{\tau}:= \frac{d}{dt}\tau$ is positive, otherwise we replace 
 $t$ by  $-t$. Then, the equation    \eqref{umparametrisation} along the geodesic reads 
  \begin{equation} 
   \phi(t) = \tfrac{1}{2} \log(\dot\tau(t)) + \const_0 \label{un2} . \end{equation}

 Now let us consider the equation \eqref{phi0}.  Substituting    
 \begin{equation} 
   \phi(t)= -\tfrac{1}{2}\log(p(t)) + \const_0  \label{un5}  \end{equation} in it (since $\dot\tau>0$, the substitution is global),   we obtain
 \begin{equation} \dddot p = 4  K g(\dot\gamma, \dot\gamma) \dot p \label{un3} . \end{equation}  Since the length of the tangent vector is preserved along a geodesic,  $g(\dot\gamma, \dot\gamma)$, and therefore $4  K g(\dot\gamma, \dot\gamma)  $ is a constant.
 The assumptions above imply that this constant is nonnegative. 
 
 Indeed, if the metric is essentially  pseudo-Riemannian, this constant is zero, since $\gamma$ is an  light-like geodesic.  If the metric is Riemannian of nonpositive    curvature,
  $g(\dot\gamma, \dot\gamma)\ge 0$, and $K\ge 0$, so their product is nonnegative.

  The equation  \eqref{un3} can be solved. 
   We will first consider the case  $Kg(\dot\gamma, \dot\gamma)= 0$- In this case, 
  the solution of \eqref{un3} is $p(t) =C_2 t^2 + C_1 t +C_0$.  Combining \eqref{un5} with \eqref{un2}, we see that  $\dot\tau = \frac{1}{C_2 t^2 + C_1 t +C_0}$.  Then
   \begin{equation} 
   \tau(t) = \int_{t_0}^t \frac{d\xi }{C_2 \xi^2 + C_1 \xi +C_0}\  \ +\const.
  \end{equation} 
  We see that if the polynomial $ C_2 t^2 + C_1 t + C_0$  has real roots (which is always the case if $C_2=0$, $C_1\ne 0$), then  the   integral 
   explodes in finite time. If the polynomial has no real roots, but $C_2\ne 0$,  the function $\tau$ is bounded. Thus, the only possibility for  $\tau $ to  be a diffeomorphism is $C_2=C_3=0$ implying  $\tau(t)  = \tfrac{1}{C_0} t  + \const_1$ implying $\dot\tau =\tfrac{1}{ C_0}$ implying $\phi$ is  constant along the  geodesic.

   Now, let us consider the case $Kg(\dot\gamma, \dot\gamma)> 0$. 
 In this case,  the general solution of the equation \eqref{un3} is 
 \begin{equation}C + C_+e^{2\sqrt{K g(\dot\gamma, \dot\gamma)}t}+ C_-e^{-2\sqrt{K g(\dot\gamma, \dot\gamma)} t}\label{un10} . \end{equation}
  Then, the function $\tau$ satisfies the ODE $\dot\tau = \frac{1}{C + C_+e^{2\sqrt{K g(\dot\gamma, \dot\gamma)}t}+ C_-e^{-2\sqrt{K g(\dot\gamma, \dot\gamma)}t}} $  implying 
  
  \begin{equation} \tau(t)  =  \int_{t_0}^t \frac{d\xi}{C + C_+e^{2\sqrt{K g(\dot\gamma, \dot\gamma)}\xi }+ C_-e^{-2\sqrt{K g(\dot\gamma, \dot\gamma)}\xi}} \ \  + \const . \label{un7} \end{equation}

  If one of the constants $C_+, C_-$ is  not zero,  the integral \eqref{un7}  is  bounded from one side, or explodes in finite time. Thus, the only possibility for $\tau $ to be a diffeomorphism of 
    $\mathbb{R}$ on itself is  $C_+=C_-=0$. Finally, $\phi$ is a constant along the geodesic $\gamma$.

    Since every point of a connected manifold 
     can be reached by a sequence of light like geodesics in the pseudo-Riemannian case,  or by a sequence of  geodesics in the Riemannian case, $\phi$ is a constant, so that $\phi_i\equiv 0$, 
     and the metrics are affine equivalent by \eqref{c1}, \qed \label{proof1} 
  
  \begin{Rem} 
  Similar idea was used by Couty  \cite{couty} in investigation of projective transformations of Einstein manifolds, and by Shen \cite{shen} in investigation of Finsler Einstein geodesically equivalent metrics.    
  \end{Rem} 
  
  \subsection{Proof of Theorem  \ref{einstein1} for Riemannian metrics of positive 
   scalar curvature} \label{proof2}

  We assume that  $g$ is a complete Einstein Riemannian 
  metric of positive scalar  curvature on a connected manifold (we do  not need that the second metric is complete).
   Then, by Corollary \ref{Tanno}, $\lambda$ is a solution of \eqref{tanno}. If the metrics are not  affine equivalent, $\lambda$ is not identically constant. 
   
   The equation  \eqref{tanno} was studied by Obata and Tanno in \cite{Obata,Tanno} in a completely different geometrical context. They proved (actually, Tanno \cite{Tanno}, because the proof of Obata \cite{Obata} has a mistake) that a complete Riemannian $g$  such that there exists a nonconstant function $\lambda $ satisfying \eqref{tanno} must have a constant positive curvature.  Applying this  result in our situation, we obtain the claim, \qed

   \section{Proof of Theorem  2\weg{\ref{dim4}}} 
  
   It is sufficient to prove Theorem \ref{dim4} in a neighborhoods of points $p$ such that $\lambda_i$ given by \eqref{lambda} does not vanish.
    Indeed,  by Corollary \ref{100}, either such points are everywhere dense, or  the metrics are affine equivalent.  We will first formulate two
     simple lemmas from Linear Algebra, then prove a simple Lemma \ref{last2}  which generalizes certain result of Levi-Civita \cite{Levi-Civita}, and then obtain Theorem \ref{dim4} as an easy corollary. 
   \subsection{Two simple  lemmas from Linear Algebra }

  We say that the vector $v^i$ lies in \emph{kernel} of the tensor $Z_{ijkl}$, if $v^i Z_{ijkl}=0$. 
  
  \begin{Lemma} Assume the tensor $Z_{ijkl}$ on $\mathbb{R}^4$ 
  has the following  symmetries:\label{la1} 
  \begin{equation} \label{sym} 
  Z_{ijkl}=Z_{klij} \ , \ Z_{ijkl}= -Z_{jikl}, 
   \end{equation}
   and satisfies $Z_{ijkl}g^{ik}=0$.   
  Suppose  the vector  $v^i$ such that $g(v,v):= v^i v^j g_{ij}\ne 0$ 
  lies in the kernel of   $Z_{ij kl}$.   Then, $Z=0$. 
  \end{Lemma} 
  \begin{Rem} The assumption    $g(v,v)\ne0 $ is important: one immediately constructs a counterexample. The dimension is also important: the claim fails for dimensions $\ge 5$. 
  \end{Rem} 
  { Proof of Lemma \ref{la1} is an easy exercise  and will be left to the reader.} We recommend to  
  consider a basis such that 
the   first vector is  $v$ and the metric is given by the matrix
  $$
  \begin{pmatrix} 
  \varepsilon_1&&& \\ &\varepsilon_2&&\\ &&\varepsilon_3&\\ &&&\varepsilon_4
  \end{pmatrix},  
  $$
   where all $\varepsilon_i\ne 0$.
   Then, the conditions $v^i Z_{ijkl}=0$ and $Z_{ijkl}g^{ik}=0$ are a system of homogeneous  linear equations on  the components of $Z$ which admits only trivial solution implying the claim, \qed

   \begin{Lemma}  \label{last} Let $a$ and $Z$ be $n\times n$ matrices    over $\mathbb{C}$ such that  
   $Z$  is skew-symmetric and such that their product $aZ$ is symmetric. Let the geometric multiplicity of the  eigenvalue 
    $\rho\in \mathbb{C}$ of the matrix $a$  be 1.   
      Then,  every vector $v$ from the generalized eigenspace  of $\rho$ lies in the kernel of the matrix $Z$.  
   \end{Lemma} 
   
   (Recall that \emph{geometric multiplicity} of $\rho$ is the dimension of the kernel of $(a-\rho \cdot {\bf 1})$ , and the  
   generalized eigenspace  of $\rho$ is the kernel of $(a- \rho \cdot {\bf 1})^n$.) 
   
{The proof of Lemma \ref{last} }  is an easy exercise in linear algebra and  will be  left to
 the reader. We recommend to consider the basis  such that the matrix
 $a$ is in   Jordan form, and then to  calculate the matrix  $aZ$. One immediately sees that it is  block diagonal,   and that if    the eigenspace  is one dimensional  then the corresponding block is  trivial,  \qed
   
   \begin{Cor} 
   Suppose $Z_{ijkl}$ is skew-symmetric with respect to indexes $i,j$. Suppose 
    \begin{equation} a_i^\alpha Z_{\alpha jkl} +a_j^\alpha Z_{\alpha i kl} =0 \label{Z} \end{equation}   for a (1,1)-tensor $a$ satisfying 
    $a_j^\alpha g_{\alpha i} = a_i^\alpha g_{\alpha j},$ where (the metric)
     $g$ is a symmetric  nondegenerate $(0,2)$-tensor. 
    We assume that all  components of $Z$, $g $, and $a$ are real.   Suppose there exists a (possible, complex) eigenvalue $\rho$ with  geometric multiplicity $1$.
     Then, there exists a vector $v$ such that $g_{ij} v^i v^j\ne 0$ lieing in the kernel of $Z$.   
   \label{cr} \end{Cor} 
   {\bf Proof. }  The condition $a_j^\alpha Z_{\alpha j} +a_i^\alpha Z_{\alpha i} =0 $ precisely means  that the matrix $aZ$ is symmetric. We see that this  condition   is the condition \eqref{Z} with ``forgotten" indexes $k$ and $l$. Then, by  Lemma  \ref{last}, every vector $v $  from the  sum of the generalized eigenspaces of  $ \rho$ and of its complex-conjugate $\bar \rho$ lies in kernel of $Z$. Since the generalized eigenspaces of $ \rho$ and of  $\bar \rho$ are orthogonal to all other generalized eigenspaces because of the condition  $a_j^\alpha g_{\alpha i} = a_i^\alpha g_{\alpha j},$ and  because the direct sum of all all generalized eigenspaces coincides with the whole vector space, the  sum of the generalized eigenspaces of  $ \rho$ and of $\bar \rho$   contains a (real) vector  $v$ such that $g_{ij} v^i v^j\ne 0$, \qed 
  
 \subsection{ If all  eigenspaces  are  more than one-dimensional, the metrics are affine equivalent.}

   \begin{Lemma}  \label{last2} 
 If geometric multiplicity of every eigenvalue of  the solution  $a_{ij}$ of \eqref{basic}   is at least two,   then the  function $\lambda$ given by \eqref{lam} is constant. 
    \end{Lemma} 
     
    \begin{Rem}
    For Riemannian metrics, the statement is due to Levi-Civita \cite{Levi-Civita}. The proof for the pseudo-Riemannian case is essentially the same, the additional difficulties are  due to possible   Jordan blocks. In a certain from, it appears in \cite{Aminova}.  
   \end{Rem}
      
    {\bf Proof of Lemma \ref{last2}.}   We prove the lemma assuming every  Jordan-Block of $a_j^i$  is as most 3-dimensional, this is sufficient for our four-dimensional  goals. The proof for arbitrary dimensions of  Jordan blocks  can be done by induction.

       Let $\rho $ be 
     an eigenvalue  of $a_{j}^i $; let $u^i$ be an eigenvector corresponding to $\rho$. 
        In a small neighborhood of almost every point,  $\rho$  a smooth (possibly, complex-valued) function. We will show that the differential $\rho_{,i}$ is proportional to $u_i$. If the eigenspace of $\rho$ is more than one-dimensional, this will imply that   $\rho_{,i}$ is constant. This implies that if all  eigenspaces  are  more than one-dimensional, the trace of $a_{j}^i$ is constant implying the metrics are affine equivalent. 
        
       Let $u$ be an eigenvector corresponding to $\rho$, i.e.,  
       \begin{equation} \label{1} 
       u_\alpha a^\alpha_i = \rho u_i  .
      \end{equation}   
  We   take  the covariant derivative and use \eqref{basic}. We obtain 
    \begin{equation} 
       u_{\alpha,j} a^\alpha_i + u_\alpha \lambda^\alpha g_{ij} + \lambda_i u_j   = 
       \rho_{,j} u_i+ \rho u_{i,j}  \label{2}.
      \end{equation}  
    We multiply \eqref{2} with $u^i$ and sum over $i$, to  obtain (using \eqref{1})  
    \begin{equation} 
       2 \lambda_\alpha u^\alpha u_j =  u_{\alpha} u^\alpha \rho_{,j}  \label{3}.
      \end{equation}  
    We see that if $u_{\alpha} u^\alpha\ne 0$ (which is in particular always the case when the Jordan block corresponding to $\rho$  is 1-dimensional), we are done. 
    
    Suppose the Jordan block corresponding to $\rho$  is more than 1-dimensional, i.e., there exists $v_i$ such that  
    \begin{equation} \label{4} 
       v_\alpha a^\alpha_i = \rho v_i + u_i  .
      \end{equation} 
Then, $u_i$ is automatically a light like vector: indeed, multiplying  \eqref{4} by $u^i$, summing  over $i$,   and using    \eqref{1},  we obtain 
    \begin{equation} \label{5} 
       u_\alpha u^\alpha = 0.
      \end{equation} 
    Differentiating \eqref{5}, we obtain 
     \begin{equation} \label{6} 
       u_{\alpha,i} u^\alpha = 0.
      \end{equation} 
 Substituting \eqref{5}  in  \eqref{3}, we obtain  
    $
       \lambda_{\alpha} u^\alpha = 0.
      $
       Differentiating \eqref{4} and using \eqref{basic}, we obtain 
     \begin{equation} \label{7} 
       v_{\alpha,j} a^\alpha_i + v_{\alpha} \lambda^\alpha g_{ij} + \lambda_iv_j=
        \rho_{,j}  v_i + \rho v_{i,j} + u_{i,j}  .
      \end{equation} 
    Multiplying \eqref{7} by $u^i$ and summing over $i$, we obtain 
    \begin{equation} \label{8} 
       v_{\alpha}  \lambda^\alpha u_j = v_\alpha u^\alpha \rho_{,j}.   
      \end{equation} 
    We see that if $v_\alpha u^\alpha \ne 0$, (which is in particular always the case when the Jordan block corresponding to $\rho$  is 2-dimensional), we are done. 
 
 Suppose the Jordan block corresponding to $\rho$  is precisely 3-dimensional, i.e., there exists $w_i$ such that \begin{equation} w_\alpha u^\alpha\ne 0\label{101} \end{equation}  and such that  
    \begin{equation} \label{9} 
       w_\alpha a^\alpha_i = \rho w_i + v_i  .
      \end{equation} 
    We multiply \eqref{9} with $v^i$ and sum over $i$, to obtain 
      \begin{equation} \label{10} 
       w_\alpha u^\alpha = v_\alpha v^\alpha  .
      \end{equation} 
     We multiply \eqref{9} with $u^i$ and sum over $i$, to obtain 
      \begin{equation} \label{11} 
       u_\alpha v^\alpha =  0 .
      \end{equation} 
     Differentiating \eqref{11}, we obtain 
      \begin{equation} \label{12} 
       u_{\alpha,i}  v^\alpha =-  u^\alpha  v_{\alpha,i} .
      \end{equation} 
     Moreover, combining \eqref{11} with \eqref{8}, we obtain  $
       \lambda_{\alpha} v^\alpha = 0.
      $
     Differentiating \eqref{9}, we obtain  $$
       w_{\alpha,j} a^\alpha_i + w_{\alpha} \lambda^\alpha g_{ij} + \lambda_iw_j=
        \rho_{,j}  w_i + \rho w_{i,j} + v_{i,j}  .$$
     Contracting this with $u^i$, we obtain 
    \begin{equation} \label{13} 
    w_{\alpha} \lambda^\alpha u_j = w_\alpha u^\alpha \rho_{,j}+ u^\alpha v_{\alpha,j} \stackrel{\eqref{12}}{=} w_\alpha u^\alpha \rho_{,j}- u_{\alpha,i} v^{\alpha} . 
    \end{equation} 
   We multiply \eqref{7} with $v^i$ and sum over $i$ to obtain   
    \begin{equation} \label{14} 
    v_{\alpha,j} u^\alpha = v_\alpha v^\alpha  \rho_{,j} + u_{\alpha, j} v^\alpha 
   . 
    \end{equation} Using \eqref{12}, we obtain 
    \begin{equation} \label{15} 
   2  v_{\alpha,j} u^\alpha = v_\alpha v^\alpha  \rho_{,j} \stackrel{\eqref{10}}{=} w_\alpha u^\alpha \rho_{,j}  
   . 
    \end{equation} 
    Combining \eqref{15} with \eqref{13}, we obtain $2 w_\alpha \lambda^\alpha u_j = 3 u_\alpha w^\alpha \rho_{,j}$. Combining this with \eqref{101}, we obtain that the differential $\rho_{,i}$ is proportional to the eigenvector $u_i$. If the eigenspace of $\rho$ is more that one-dimensional, this implies that $\rho_{,i}=0$,   \qed

    \subsection{Proof of Theorem {\ref{dim4}}} 
   If the dimension is 3, Theorem  \ref{dim4} follows from the well-known fact that every Einstein 3-manifold has constant curvature.

    We  assume that $g$ is an  Einstein metric on $M^4$. Let  $\bar 
    g$ be geodesically equivalent to $g$. We consider the solution $a_{ij}$ of \eqref{basic} given by \eqref{a}.  Assume that the correspondent $\lambda_i\ne 0$ at $p$.  We will show that in a small neighborhood of $p$ the metric $g$ has constant curvature implying the metrics $\bar g$ and $\hat g$ have constant curvature as well by Beltrami  Theorem  (see for example \cite{beltrami_short}, or the original papers \cite{Beltrami} and \cite{schur}).
          
 Substituting  equation \eqref{vb}  in \eqref{int1}, we obtain
$  a_{i \alpha} Z^{\alpha}_{jkl} +  a_{\alpha j}Z^{\alpha}_{ikl}=0, 
$
 where \begin{equation} \label{view} Z^{i}_{jkl}= R^i_{jkl} - K\cdot  ( \delta^i_lg_{jk} - \delta^i_kg_{jl}  ).\end{equation} We see that by construction the tensor $Z_{ijkl}$ has the  symmetries \eqref{sym}. Since $g$ is Einstein, the  tensor $Z_{ijkl}$  satisfies $Z_{ijkl}g^{ik}=0$. 
 
 By Lemma \ref{last2}, at almost every point there exists an eigenvalue  of $a_{j}^{i} $ with geometric multiplicity one.   Then, by  Corollary \ref{cr},  there exists a vector $v^i$ such that $g(v,v) \ne 0$ and such that $v^i$ lies in the kernel of $Z$. By Lemma \ref{last}, the tensor $Z\equiv 0$ implying in view of \eqref{view} the claim, \qed

\end{document}